\documentclass[11pt]{article}
\usepackage{amssymb,amsmath,latexsym, epsfig}

\begin{document}

\begin{doublespace}

\newtheorem{thm}{Theorem}[section]
\newtheorem{lemma}[thm]{Lemma}
\newtheorem{defn}{Definition}[section]
\newtheorem{prop}[thm]{Proposition}
\newtheorem{corollary}[thm]{Corollary}
\newtheorem{remark}[thm]{Remark}
\newtheorem{example}[thm]{Example}
\numberwithin{equation}{section}

\def\ee{\varepsilon}
\def\qed{{\hfill $\Box$ \bigskip}}
\def\MM{{\cal M}}
\def\BB{{\cal B}}
\def\LL{{\cal L}}
\def\FF{{\cal F}}
\def\GG{{\cal G}}
\def\EE{{\cal E}}
\def\QQ{{\cal Q}}

\def\R{{\mathbb R}}
\def\L{{\bf L}}
\def\E{{\mathbb E}}
\def\F{{\bf F}}
\def\P{{\mathbb P}}
\def\N{{\mathbb N}}
\def\eps{\varepsilon}
\def\wh{\widehat}
\def\pf{\noindent{\bf Proof.} }

\title{\Large \bf Estimation of Parameters of Stable Distributions  }
\author{ Chunlin Wang\\
Department of Mathematics\\
University of Illinois\\
Urbana, IL 61801, USA\\
Email: cwang13@uiuc.edu}
\date{}
\maketitle
\date{}
\maketitle

\begin{abstract}
In this paper, we propose a method based on GMM (the generalized
method of moments) to estimate the parameters of stable
distributions with $0<\alpha<2$. We don't assume symmetry for
stable distributions.
\end{abstract}



 \vspace{.1truein}
\noindent{\bf Running Title:} Estimation of parameters of stable
distributions

\bigskip

{ \section {{Introduction}}}

The class of stable distributions is an important probability
distribution class, which can be viewed as the sum of a large
number of independent and identically distributed random variables
with very small individual effects. The stable non-Gaussian
distributions has the property that the second moment is infinite.
Plenty of evidence suggested that some important economic
variables such as stock price changes, interest rate changes, and
price expectations etc. can be better described by stable
distributions, in most cases, by stable non-Gaussian
distributions, see Mandelbrot (1963a,b), Fama (1965) etc.. It was
suggested by Bartels (1977), Koenker and Bassett (1978) that the
distribution of the regression disturbance may also belong to the
class of stable distributions

One dimensional stable distributions with parameter
 $\alpha=1$ and $\alpha=2$ are Cauchy distribution and normal distribution
 respectively. Both are well studied.

 For stable distributions with $0<\alpha<2$, since there are no closed
 form expressions for the density function in most cases, estimation of
 parameters encounters difficulties. Some important
 works in this area were done by Fama and Roll (1968,1971) and others. Stable distributions were also used in tests for nomality by Bera and
 McKenzie (1986).

 Previous methods on estimation of parameters of stable distributions include the fractile method of Famma and Roll (1968,1971), the improved version of McCulloch
 (1986), the approximate maximum likelihood theory developed by DuMouchel (1973a,b,1975), and the iterative regression method of
 Koutrouvelis (1980,1981) etc.. Akgiray and Lamoureux (1989) made a comparative study of the
 fractile method and iterative regression method.

 In this paper we propose a method based on GMM (the generalized method of moments) to estimate simultaneously all the parameters of
  stable distributions with $1<\alpha<2$ and $0<\alpha<1$. We don't assume the
 symmetry of stable distribution here.

  Stable distributions with $0<\alpha<2$ and $\alpha\neq 1$ may be defined by the
 characteristic function
 \begin{equation}
  \hat{\mu}=\exp[-c|z|^{\alpha}(1-i\beta \textrm{tan}\frac{\pi
 \alpha}{2}\textrm{sgn}
 z)+i\tau z],
 \end{equation}
 with $c>0,\beta \in [-1,1] $ and $\tau \in \R$, where tan is tangent function and sgn is sign function.  $\beta=0$ corresponds to symmetric stable
 distributions. $\beta \neq 0$ corresponds to nonsymmetric
 stable distributions. Two special cases $\beta=1$ and $\beta=-1$ correspond to one sided stable distributions. Here we focus on the cases that  $-1<\beta<1$. $\tau$ is the drift term. For detailed
 account of stable distributions,
 see Sato (1999) and Zolotarev (1983).

We assume the random variable $x$ has the stable
 distribution with parameters $\alpha, \beta, \tau,$ and $c$.

In order to present infinite and finite series expressions of the
density function of the stable distribution in Sato (1999), we
introduce a new parameter $\tilde{c}$ which satisfies
 \begin{eqnarray*}
 &&c=\cos(\pi\beta\frac{2-\alpha}{2})\tilde{c}^{\alpha}, \,\,\text{ if $1<\alpha<2$;}\\
&&c=\cos(\frac{\pi\beta\alpha}{2})\tilde{c}^{\alpha},\,\,\,\,\,\,\,\,\,\,\,\,
\text{ if $0<\alpha<1$.}
\end{eqnarray*}

 Thus estimation of $\tilde{c}$ is equivalent to
estimation of $c$. Consider the new variable
$$u=\frac{x-\tau}{\tilde{c}}.$$ $u$ has
$\cos(\pi\beta\frac{2-\alpha}{2})$ if $1<\alpha<2$ and
$\cos(\frac{\pi\beta\alpha}{2})$ if $0<\alpha<1$  as a parameter
which appears in the position of $c$ in the characteristic
function expression (1.1).

By Sato (1999), we have the following infinite and finite series
expressions of the density function associated with variable $u$.

\begin{enumerate}
\item{When $1<\alpha<2$}, the convergent series expression for the
density function is
 \begin{equation}
p_{\alpha}(u)=\frac{1}{\pi}\sum_{n=1}^{\infty}(-1)^{n-1}\frac{\Gamma(n/\alpha+1)}{n!}(\sin(\pi
n(1-\beta\frac{2-\alpha}{\alpha})/2))u^{n-1}, \,\,\,\,\text{for
$u\in \R $,}
\end{equation}

 and for $u\rightarrow \infty$, the finite series expression is:
\begin{equation}
\begin{split}
p_\alpha(u)&=\frac{1}{\pi}\sum_{n=1}^{N}(-1)^{n-1}\frac{\Gamma(n\alpha+1)}{n!}(\sin(\pi
n(1-\beta\frac{2-\alpha}{\alpha})/2\alpha))u^{-n\alpha-1}\\
            &\textrm{ }+O(u^{-(N+1)\alpha-1}).\,\,\,\,\,
\end{split}
\end{equation}

For $u\rightarrow -\infty$, since the dual distribution
$\tilde{u}=-u$ has the density $\tilde{p}_\alpha(u)=p_\alpha(-u)$,
the finite series expression is:
\begin{equation}
\begin{split}
p_\alpha(u)&=\frac{1}{\pi}\sum_{n=1}^{N}(-1)^{n-1}\frac{\Gamma(n\alpha+1)}{n!}(\sin(\pi
n(1+\beta\frac{2-\alpha}{\alpha})/2 \alpha))|u|^{-n\alpha-1}\\
            &\textrm{ }+O(|u|^{-(N+1)\alpha-1}).\\
\end{split}
\end{equation}

\item{When $0<\alpha<1$}, for $u>0$, the convergent series
expression for the density function is
 \begin{equation}
p_{\alpha}(u)=\frac{1}{\pi}\sum_{n=1}^{\infty}(-1)^{n-1}\frac{\Gamma(n\alpha+1)}{n!}(\sin(\pi
n(1+\beta)/2\alpha))u^{-n\alpha-1},
\end{equation}
and for $u<0$, the convergent series expression for the density
function is
 \begin{equation}
p_{\alpha}(u)=\frac{1}{\pi}\sum_{n=1}^{\infty}(-1)^{n-1}\frac{\Gamma(n\alpha+1)}{n!}(\sin(\pi
n(1-\beta)/2\alpha))|u|^{-n\alpha-1},
\end{equation}
which is obtained by duality argument similar to the above.

For $u\rightarrow 0$, the finite series expression is:
\begin{equation}
p_\alpha(u)=\frac{1}{\pi}\sum_{n=1}^{N}(-1)^{n-1}\frac{\Gamma(n/\alpha+1)}{n!}(\sin(\pi
n(1+\beta)/2))u^{n-1}+O(u^{N}).\,\,\,\,\,
\end{equation}

\end{enumerate}

 In the following two sections, we discuss a method based on GMM to estimate
parameters $\alpha, \beta, \tau,$ and $\tilde{c}$. For
convenience, we denote the vector $(\alpha, \beta, \tau,
\tilde{c})$ by $\theta$ and the vector ($\alpha_0, \beta_0,
\tau_0, \tilde{c}_0)$ by $\theta_0$. In section 2, we look at the
case when $1<\alpha<2$. In section 3, we study the case when
$0<\alpha<1$.

\bigskip
\bigskip
\bigskip

 {\section {{$1<\alpha<2$}}}

Let us consider the following moment functions:
\begin{equation}
f_k(y;\theta)=
\begin{cases}
   d_1(\frac{y-\tau}{\tilde{c}})^{\frac{1}{k+1}},  & \text
   { if $y\geq R_1$;}\\
   \\
   (\frac{y-\tau}{\tilde{c}})^{2(1+\frac{1}{2(k+1)})}, & \text
   { if $R_2< y< R_1$;}\\
   \\
   d_2(\frac{\tau-y}{\tilde{c}})^{\frac{1}{(k+1)}},  & \text
   { if $y\leq R_2$,}\\
\end{cases}
\end{equation}
 for $k=1,2,...,m$ with $m > 4$, where the exponent
 $2(1+\frac{1}{2(k+1)})$ means first take square and then take
 exponent $1+\frac{1}{2(k+1)}$. $R_1$ and $R_2$ are
constants such that $R_2< \tau < R_1$ with
$\frac{R_1-\tau}{\tilde{c}}$ and $\frac{\tau-R_2}{\tilde{c}}$ not
large. The reason for this requirement is that the convergent
formula (1.2) of the density function is appropriate for small $u$
and the finite series expressions (1.3) and (1.4) are appropriate
for $u$ is large. This phenomenon was pointed out by Famma and
Roll (1968) for symmetric stable distribution. When we solve the
optimization problem resulted by GMM, if the solution for $\tau$
is very close to $R_1$ or $R_2$, it means that we may not get the
globally optimization yet. We have to make $R_1$ bigger or $R_2$
bigger, and do optimization again with the solution from previous
step as the initial point. About $d_1$ and $d_2$ in the above
expression, they are constants satisfying
$$d_1=((R_1-\tau)/\tilde{c})^2\,\,\text{ and }\,\, d_2=((\tau-R_2)/\tilde{c})^2.$$
We can see that $f_k(y;\theta),\,k=1,...,m,$ are continuous
functions.

Let $x$ be the random variable of stable distribution with
parameters $\theta$. Then
\begin{equation}
E(f_k(x;\theta))=\int_{-\infty}^{\infty}f_k(y;\theta)p_\alpha(\frac{y-\tau}{\tilde{c}})\frac{1}{\tilde{c}}dy,
\end{equation}
where $p_\alpha$ is the density function in (1.1), (1.2) and
(1.3).

Let $x_1,x_2,...,x_n$ be independent stable distributions with the
same parameter vector $\theta_0$. Since
$E|f_k(x_i;\theta)|<\infty,\,\,i=1,...,n,\,k=1,...,m$, by strong
law of large numbers, we have
$$\frac{\sum_{i=1}^{n}f_k(x_i;\theta_0)}{n}\rightarrow
E(f_k(x;\theta_0))\,\,\text{almost surely},\,\,k=1,...,m,$$ as
$n\rightarrow \infty$.

 Denote the sample moment
$\frac{\sum_i^nf_k(x_i;\theta))}{n}$ by
$E_nf_k(x;\theta),\,k=1,...,m.$

We can see that $E_nf_k(x;\theta))$ and $Ef_k(x;\theta)$ have
first derivatives (vector) with respect to $\theta=(\alpha, \beta,
\tau, \tilde{c})$.  Define
\begin{equation}
V=c_1[E_nf_1(x;\theta_0)-Ef_1(x;\theta))]^2+\cdots+c_m[E_nf_m(x;\theta_0)-Ef_m(x;\theta)]^2,
\end{equation} where $c_1,...,c_m$ are positive constants. Then by GMM,
the minimization solution $\hat{\theta}$ of $V$ converges in
distribution to $\theta_0$.

Since we don't have closed form expression for the density of
stable distributions, we need approximate finite series
expressions for $Ef_k(x;\theta)),\,k=1,...,m$. When we minimize
the function $V$ with respect to the parameter vector $\theta$, we
use the approximate finite series expressions for $Ef_k(x;\theta)$
instead of $Ef_k(x;\theta),\,k=1,...,m$ themselves.

In the following, we will show the existence of approximate finite
series expressions for $Ef_k(x;\theta)),\,\,k=1,...,m$.

For positive integer $N>0$, define $p_{1,\,N},$ $ p_{2,\,N}$ and
$p_{3,\,N}$ as follows:
\begin{equation}
p_{1,\,N}(u)=\frac{1}{\pi}\sum_{n=1}^{N}(-1)^{n-1}\frac{\Gamma(n\alpha+1)}{n!}(\sin(\pi
n(1-\beta\frac{2-\alpha}{\alpha})/2\alpha))u^{-n\alpha-1,}\,\,\,\,\,\,\textrm{
for }u>0,
\end{equation}
\begin{equation}
p_{2,\,N}(u)=\frac{1}{\pi}\sum_{n=1}^{N}(-1)^{n-1}\frac{\Gamma(n\alpha+1)}{n!}(\sin(\pi
n(1+\beta\frac{2-\alpha}{\alpha})/2
\alpha))|u|^{-n\alpha-1},\,\,\,\textrm{ for }u<0,
\end{equation}
\begin{equation}
p_{3,\,N}(u)=\frac{1}{\pi}\sum_{n=1}^{N}(-1)^{n-1}\frac{\Gamma(n/\alpha+1)}{n!}(\sin(\pi
n(1-\beta\frac{2-\alpha}{\alpha})/2))u^{n-1},\,\,\,\,\,\,\,\,\,\,\,\,\textrm{
for }u\in \R,
\end{equation}
with $p_{1,\,N},$ $ p_{2,\,N}$ and $p_{3,\,N}$ corresponding to
(1.3), (1.4) and (1.2) respectively.

By (1.3), we have
\begin{equation}
\begin{split}
&\int_{R_1}^{\infty}f_k(y;\theta)|p_\alpha((y-\tau)/\tilde{c})-p_{1,\,N}((y-\tau)/\tilde{c})|\frac{1}{\tilde{c}}dy\\
&\leq
M_1\int_{\frac{R_1-\tau}{\tilde{c}}}^{\infty}u^{\frac{1}{k+1}}u^{-(N+1)\alpha-1}du\,\,\,\,\,\,\,\\
&\,\,\,\,\,\,(\text{ where $M_1$ is a constant })\\
&\leq
M_1\int_{\frac{R_1-\tau}{\tilde{c}}}^{\infty}u^{-(N+1)\alpha-1+\frac{1}{k+1}}du\,\,\,\,\,\\
&\leq \text{
O}(((R_1-\tau)/\tilde{c})^{-(N+1)\alpha+\frac{1}{k+1}}),
\end{split}
\end{equation}
which goes to $0$ as $N\rightarrow \infty,$ when
$(R_1-\tau)/\tilde{c}>1.$

Similarly, by (1.4), we have
\begin{equation}
\begin{split}
&\int_{-\infty}^{R_2}f_k(y;\theta)|p_\alpha((y-\tau)/\tilde{c})-p_{2,\,N}((y-\tau)/\tilde{c})|\frac{1}{\tilde{c}}dy\\
&\leq \text{
O}(((\tau-R_2)/\tilde{c})^{-(N+1)\alpha+\frac{1}{k+1}}),
\end{split}
\end{equation}
which goes to $0$ as $N\rightarrow \infty$ when
$(\tau-R_2)/\tilde{c}>1.$

By (1.2), we have
$$p_{1,\,N}((y-\tau)/\tilde{c})\rightarrow
p_\alpha((y-\tau)/\tilde{c})\,\,\,\,\text{pointwisely for $y\in
[R_2, R_1]$,\,\, as $ N\rightarrow \infty$},$$ and the convergence
is uniform when both $((R_1-\tau)/\tilde{c})$ and
$((\tau-R_2)/\tilde{c})$ are less than $1$.

In fact, we can show that for given $\theta=(\alpha, \beta, \tau,
\tilde{c})$, given $R_1$ and $R_1$, without the assumption that
$((R_1-\tau)/\tilde{c})$ and $((\tau-R_2)/\tilde{c})$ are less
than $1$, \begin{equation}
p_{1,\,N}((y-\tau)/\tilde{c})\rightarrow
p_\alpha((y-\tau)/\tilde{c})\,\,\,\,\text{uniformly for all $y\in
[R_2, R_1]$, \,\, as $ N\rightarrow \infty$   }\end{equation}

To show this statement, we need the following lemma:
\begin{lemma}
for given $M>0$, $L>1$, there exists a positive integer $n_0$,
which depends on $L$, $M$ and $1<\alpha<2$, such that when $n\geq
n_0$,
\begin{equation}
\frac{\Gamma(n/\alpha+1)}{n!}M^{n-1}\leq L^{-n}.
\end{equation}
\end{lemma}
\pf We know that \begin{equation} \Gamma(n/\alpha+1)\leq
C_0([n/\alpha]+1)!, \end{equation} where $[n/\alpha]$ is the
greatest integer less than $n/\alpha$ and $C_0$ is a positive
constant.

Since \begin{equation} \frac{n!}{m!(n-m)!}\geq 1,\,\,\text{ for
}m=0,1,...,n, \end{equation} letting $m=[n/\alpha]+1$, then
\begin{equation}
\frac{([n/\alpha]+1)!}{n!}\leq
\frac{1}{(n-[n/\alpha]-1)!},
\end{equation}
which means
\begin{equation}
\frac{\Gamma(n/\alpha+1)}{n!}M^{n-1}\leq
C_0\frac{([n/\alpha]+1)!}{n!}M^{n-1}\leq
C_0\frac{1}{(n-[n/\alpha]-1)!}M^{n-1}.
\end{equation}
By Stirling's formula, \begin{equation} (n-[n/\alpha]-1)!\thicksim
(n-[n/\alpha]-1)^{(n-[n/\alpha]-1)}e^{-(n-[n/\alpha]-1)}\sqrt{2\pi(n-[n/\alpha]-1)},
\end{equation}
as $n\rightarrow \infty$.

Thus, in order to show the lemma, it is sufficient to show that
there exists $n_0>0$ such that for any $n\geq n_0$,
\begin{equation}
C_0M^{n-1}L^n<
(n-[n/\alpha]-1)^{(n-[n/\alpha]-1)}e^{-(n-[n/\alpha]-1)}\sqrt{(n-[n/\alpha]-1)}.
\end{equation}

It is clear that
\begin{equation}
\begin{split}
&\,\,\,\,\,\,\,\log\{(n-[n/\alpha]-1)^{(n-[n/\alpha]-1)}e^{-(n-[n/\alpha]-1)}\sqrt{2\pi(n-[n/\alpha]-1)}\}\\
&=(n-[n/\alpha]-1)\log(n-[n/\alpha]-1)-(n-[n/\alpha]-1)+\frac{1}{2}\log(2\pi)\\
&\,\,\,+\frac{1}{2}\log{(n-[n/\alpha]-1)},
\end{split}
\end{equation}
where $(n-[n/\alpha]-1)\log(n-[n/\alpha]-1)$ becomes the dominant
term when $n$ is large enough.

We can see that when $n$ is large enough,
\begin{equation}
\log{C_0}+(n-1)\log{M}+n\log{L}<(n-[n/\alpha]-1)\log(n-[n/\alpha]-1),
\end{equation}
which means that there exists $n_0>0$ such that for any $n\geq
n_0$, (2.16) holds. Therefore we approved the lemma. \qed

Let $M=\min((R_1-\tau)/\tilde{c},(\tau-R_2)/\tilde{c})$.
Substituting $\frac{y-\tau}{\tilde{c}}$ for $u$ in (1.2), applying
Lemma 2.1 and using the following fact,
$$\sum_{n=N}^{\infty}L^{-n}=
\frac{L^{-N}}{1-\frac{1}{L}},$$ we have (2.9) holds.

One thing we need to pay attention is that although we have
uniform convergence in (2.9), from the proof of Lemma 1, we can
see that the convergence speed depends on $(R_1-\tau)/\tilde{c},
(\tau-R_2)/\tilde{c}$ and $\alpha$. When $\alpha$ is close to $1$
or either $R_1$ or $R_2$ is very large, the convergence speed will
be slow. Since we cannot change the parameter $\alpha$, $\tau$ and
$\tilde{c}$, the only way to obtain good convergence speed is to
keep $R_1$ and $R_2$ close to $\tau$.

Clearly (2.9) implies that
\begin{equation}
\int_{R_2}^{R_1}f_k(y;\theta)|p_\alpha((y-\tau)/\tilde{c})-p_{3,\,N}((y-\tau)/\tilde{c})|\frac{1}{\tilde{c}}dy\rightarrow
0,
\end{equation}
as $N\rightarrow \infty.$

Define
\begin{equation}
\begin{split}
&T_{k,1,N}=\int_{R_1}^{\infty}f_k(y;\theta)p_{1,\,N}((y-\tau)/\tilde{c})\frac{1}{\tilde{c}}dy,\\
&T_{k,2,N}=\int_{-\infty}^{R_2}f_k(y;\theta)p_{2,\,N}((y-\tau)/\tilde{c})\frac{1}{\tilde{c}}dy,\\
&T_{k,3,N}=\int_{R_2}^{R_1}f_k(y;\theta)p_{3,\,N}((y-\tau)/\tilde{c})\frac{1}{\tilde{c}}dy,\\
\end{split}
\end{equation}
for $k=1,...,m$.

Combining $(2.7$, $(2.8)$ and $(2.19)$, we have
\begin{equation}
T_{k,1,N}+T_{k,2,N}+T_{k,3,N}\rightarrow E(f_k(x;\theta)),
\end{equation}
for $k=1,...,m$, as $N$ goes to $\infty$. This means the sum of
$T_{k,1,N}$, $T_{k,2,N}$, and $T_{k,3,N}$ is the approximate
finite series expression we wanted.

Next we give the expressions of $T_{k,1,N}$, $T_{k,2,N}$, and
$T_{k,3,N}$ without integrals inside. It is easy to see
\begin{enumerate}
\item \begin{equation}
\begin{split}
&\,\,\,\,\,\,\,\,\,\,T_{k,1,N}\\
&=\int_{R_1}^{\infty}f_k(y;\theta)p_{1,\,N}((y-\tau)/\tilde{c})\frac{1}{\tilde{c}}dy\\
&=\frac{(R_1-\tau)^2}{{\tilde{c}}^2\pi}\sum_{n=1}^{N}(-1)^{n-1}\frac{\Gamma(n\alpha+1)}{n!}(\sin(\pi
n(1-\beta\frac{2-\alpha}{\alpha})/2\alpha))\\
&\,\,\,\,\,\,\,\cdot \int_{\frac{R_1-\tau}{\tilde{c}}}^{\infty}u^{\frac{1}{k+1}}u^{-n\alpha-1}du,\\
\end{split}
\end{equation}
where
\begin{equation}
\begin{split}
&\,\,\,\,\,\,\,\,\int_{\frac{R_1-\tau}{\tilde{c}}}^{\infty}u^{\frac{1}{k+1}}u^{-n\alpha-1}du\\
&=\int_{\frac{R_1-\tau}{\tilde{c}}}^{\infty}u^{-n\alpha-1+\frac{1}{k+1}}du\\
&=\frac{1}{-n\alpha+\frac{1}{k+1}}((R_1-\tau)/\tilde{c})^{-n\alpha+\frac{1}{k+1}}.\\
\end{split}
\end{equation}

\item \begin{equation}
\begin{split}
&\,\,\,\,\,\,\,\,\,\,T_{k,2,N}\\
&=\int_{-\infty}^{R_2}f_k(y;\theta)p_{2,\,N}((y-\tau)/\tilde{c})\frac{1}{\tilde{c}}dy\\
&=\frac{((\tau-R_2)/\tilde{c})^{2}}{\pi}\sum_{n=1}^{N}(-1)^{n-1}\frac{\Gamma(n\alpha+1)}{n!}(\sin(\pi
n(1+\beta\frac{2-\alpha}{\alpha})/2\alpha))\\
&\,\,\,\,\,\,\,\cdot \int_{-\infty}^{\frac{R_2-\tau}{\tilde{c}}}|u|^{\frac{1}{k+1}}|u|^{-n\alpha-1}du,\\
\end{split}
\end{equation}
where
\begin{equation}
\begin{split}
&\,\,\,\,\,\,\,\,\int_{-\infty}^{\frac{R_2-\tau}{\tilde{c}}}|u|^{\frac{1}{k+1}}|u|^{-n\alpha-1}du\\
&=\int_{\frac{\tau-R_2}{\tilde{c}}}^{\infty}u^{-(n+1)+\frac{1}{k+1}}du\\
&=\frac{1}{-n\alpha+\frac{1}{k+1}}((\tau-R_2)/\tilde{c})^{-n\alpha+\frac{1}{k+1}}.\\
\end{split}
\end{equation}

\item \begin{equation}
\begin{split}
&\,\,\,\,\,\,\,\,\,\,T_{k,3,N}\\
&=\int_{R_1}^{R_2}f_k(y;\theta)p_{3,\,N}((y-\tau)/\tilde{c})\frac{1}{\tilde{c}}dy\\
&=\frac{1}{\pi}\sum_{n=1}^{N}(-1)^{n-1}\frac{\Gamma(n/\alpha+1)}{n!}(\sin(\pi
n(1-\beta\frac{2-\alpha}{\alpha})/2))\\
&\,\,\,\,\,\,\,\cdot \int_{\frac{R_2-\tau}{\tilde{c}}}^{\frac{R_1-\tau}{\tilde{c}}}u^{2(1+\frac{1}{2(k+1)})}u^{n-1}du,\\
\end{split}
\end{equation}
\end{enumerate}
where
\begin{equation}
\int_{\frac{R_2-\tau}{\tilde{c}}}^{\frac{R_1-\tau}{\tilde{c}}}u^{2(1+\frac{1}{2(k+1)})}u^{n-1}du=\\
\begin{cases}
A_1+A_2, &\text{ if $n$ is odd;}\\
A_1-A_2, &\text{ if $n$ is even,}\\
\end{cases}
\end{equation}
with
\begin{equation}
\begin{split}
A_1&=\int_{0}^{\frac{R_1-\tau}{\tilde{c}}}u^{2(1+\frac{1}{2(k+1)})}u^{n-1}du\\
   &=\int_{0}^{\frac{R_1-\tau}{\tilde{c}}}u^{n-1+2(1+\frac{1}{2(k+1)})}du\\
   &=\frac{1}{n+2(1+\frac{1}{2(k+1)})}((R_1-\tau)/\tilde{c})^{n+2(1+\frac{1}{2(k+1)})},
\end{split}
\end{equation}
and
\begin{equation}
\begin{split}
A_2&=\int_{0}^{\frac{\tau-R_2}{\tilde{c}}}u^{2(1+\frac{1}{2(k+1)})}u^{n-1}du\\
   &=\int_{0}^{\frac{\tau-R_2}{\tilde{c}}}u^{n-1+2(1+\frac{1}{2(k+1)})}du\\
   &=\frac{1}{n+2(1+\frac{1}{2(k+1)})}((\tau-R_2)/\tilde{c})^{n+2(1+\frac{1}{2(k+1)})}.
\end{split}
\end{equation}

Thus we obtain the expressions of $T_{k,1,N}$, $T_{k,2,N}$, and
$T_{k,3,N}$ without integrals inside. By (2.21), we can substitute
the sum of $T_{k,1,N}$, $T_{k,2,N}$ and $T_{k,3,N}$ for
$E(f_k(x;\theta))$ in the expression of $V$ of (2.3) and then do
optimization to get the estimate $\hat{\theta}$ for $\theta_0$.

\bigskip
\bigskip
\bigskip

{\section {$0<\alpha<1$}}

In order to ensure integrability, we consider the following moment
functions, which are different from those when $1<\alpha<2$:
\begin{equation}
\bar{f}_k(y;\theta)=
\begin{cases}
   \bar{d}_1(\frac{y-\tau}{\tilde{c}})^{\frac{-1}{k+1}},  & \text
   { if $y\geq \bar{R}_1$;}\\
   \\
   (\frac{y-\tau}{\tilde{c}})^{2(1-\frac{1}{2(k+1)})}, & \text
   { if $\bar{R}_2 < y < \bar{R}_1$;}\\
   \\
   \bar{d}_2(\frac{\tau-y}{\tilde{c}})^{\frac{-1}{(k+1)}},  & \text
   { if $y\leq \bar{R}_2$,}\\
\end{cases}
\end{equation}
 for $k=1,2,...,m$ with $m > 4$, where the exponent
 $2(1-\frac{1}{2(k+1)})$ means first take square and then take
 exponent $1-\frac{1}{2(k+1)}$. $\bar{R}_1$ and $\bar{R}_2$ are
constants such that $\bar{R}_2< \tau < \bar{R}_1$ with
$\frac{\bar{R}_1-\tau}{\tilde{c}}$ and
$\frac{\tau-\bar{R}_2}{\tilde{c}}$ are small. The reason for this
requirement is that when $u$ is small, we have the finite series
expression (1.7) for the density function. $\bar{d}_1$ and
$\bar{d}_2$ in the above expression are constants satisfying
$$\bar{d}_1=((\bar{R}_1-\tau)/\tilde{c})^2\,\,\text{ and }\,\, \bar{d}_2=((\tau-\bar{R}_2)/\tilde{c})^2.$$
It is clear that $\bar{f}_k(y;\theta),\,k=1,...,m,$ are continuous
functions.

Let $x$ be the random variable of stable distribution with
parameters $\theta$. Then
\begin{equation}
E(\bar{f}_k(x;\theta))=\int_{-\infty}^{\infty}\bar{f}_k(y;\theta)p_\alpha(\frac{y-\tau}{\tilde{c}})\frac{1}{\tilde{c}}dy,
\end{equation}
where $p_\alpha$ is the density function in (1.5), (1.6) and
(1.7).

Similar to the argument in section 2, for independent stable
distributions  $x_1,x_2,...,x_n$ with the same parameter vector
$\theta_0$, denote the sample moment
$\frac{\sum_i^n\bar{f}_k(x_i;\theta))}{n}$ by
$E_n\bar{f}_k(x;\theta),\,k=1,...,m.$ Since
$E|\bar{f}_k(x_i;\theta)|<\infty,\,\,i=1,...,n,\,k=1,...,m$, by
strong law of large numbers, we have
$$E_n\bar{f}_k(x;\theta_0)\rightarrow
E(\bar{f}_k(x;\theta_0))\,\,\text{almost surely},\,\,k=1,...,m,$$
as $n\rightarrow \infty$.

Both $E_n\bar{f}_k(x;\theta))$ and $E\bar{f}_k(x;\theta)$ have
first derivatives (vector) with respect to $\theta=(\alpha, \beta,
\tau, \tilde{c})$.  We define
\begin{equation}
\bar{V}=\bar{c}_1[E_n\bar{f}_1(x;\theta_0)-E\bar{f}_1(x;\theta))]^2+\cdots+\bar{c}_m[E_n\bar{f}_m(x;\theta_0)-E\bar{f}_m(x;\theta)]^2,
\end{equation} where $\bar{c}_1,...,\bar{c}_m$ are positive constants. Then by GMM,
the minimization solution $\hat{\theta}$ of $\bar{V}$ converges in
distribution to $\theta_0$.

Because of no closed form expressions for the density of stable
distributions, we substitute its approximate finite series
expression (if they exist) for $E\bar{f}_k(x;\theta)),\,k=1,...,m$
in the expression of $\bar{V}$ of (3.3).

Using the similar argument as we did in section 2, it is not hard
to show the existence of approximate finite series expression for
$E\bar{f}_k(x;\theta)),\,\,k=1,...,m$. For $N>0$, we denote it by
the sum of $\bar{T}_{k,1,N}$, $\bar{T}_{k,2,N}$, and
$\bar{T}_{k,3,N}$, where
\begin{equation}
\begin{split}
&\bar{T}_{k,1,N}=\int_{\bar{R}_1}^{\infty}\bar{f}_k(y;\theta)\bar{p}_{1,\,N}((y-\tau)/\tilde{c})\frac{1}{\tilde{c}}dy,\\
&\bar{T}_{k,2,N}=\int_{-\infty}^{\bar{R}_2}\bar{f}_k(y;\theta)\bar{p}_{2,\,N}((y-\tau)/\tilde{c})\frac{1}{\tilde{c}}dy,\\
&\bar{T}_{k,3,N}=\int_{\bar{R}_2}^{\bar{R}_1}\bar{f}_k(y;\theta)\bar{p}_{3,\,N}((y-\tau)/\tilde{c})\frac{1}{\tilde{c}}dy,\\
\end{split}
\end{equation}
for $k=1,...,m$, with $\bar{p}_{1,\,N},$ $ \bar{p}_{2,\,N}$ and
$\bar{p}_{3,\,N}$ being:
\begin{equation}
\bar{p}_{1,\,N}(u)=\frac{1}{\pi}\sum_{n=1}^{N}(-1)^{n-1}\frac{\Gamma(n\alpha+1)}{n!}(\sin(\pi
n(1+\beta)/2\alpha))u^{-n\alpha-1},     \,\,\,\,\,\,\textrm{ for
}u>0,
\end{equation}
\begin{equation}
\bar{p}_{2,\,N}(u)=\frac{1}{\pi}\sum_{n=1}^{N}(-1)^{n-1}\frac{\Gamma(n\alpha+1)}{n!}(\sin(\pi
n(1-\beta)/2\alpha))|u|^{-n\alpha-1}   ,\,\,\,\textrm{ for }u<0,
\end{equation}
\begin{equation}
\bar{p}_{3,\,N}(u)=
\frac{1}{\pi}\sum_{n=1}^{N}(-1)^{n-1}\frac{\Gamma(n/\alpha+1)}{n!}(\sin(\pi
n(1+\beta)/2))u^{n-1}    ,\,\,\,\,\,\,\,\,\,\,\,\,\textrm{ for
}u\in \R,
\end{equation}
We can see that $\bar{p}_{1,\,N},$ $ \bar{p}_{2,\,N}$ and
$\bar{p}_{3,\,N}$ corresponding to (1.5), (1.6) and (1.7)
respectively.

Next we give the expressions of $\bar{T}_{k,1,N}$,
$\bar{T}_{k,2,N}$, and $\bar{T}_{k,3,N}$ without integrals inside.
It is easy to see
\begin{enumerate}
\item \begin{equation}
\begin{split}
&\,\,\,\,\,\,\,\,\,\,\bar{T}_{k,1,N}\\
&=\int_{\bar{R}_1}^{\infty}\bar{f}_k(y;\theta)\bar{p}_{1,\,N}((y-\tau)/\tilde{c})\frac{1}{\tilde{c}}dy\\
&=\frac{(\bar{R}_1-\tau)^2}{{\tilde{c}}^2\pi}\frac{1}{\pi}\sum_{n=1}^{N}(-1)^{n-1}\frac{\Gamma(n\alpha+1)}{n!}(\sin(\pi
n(1+\beta)/2\alpha))      \\
&\,\,\,\,\,\,\,\cdot \int_{\frac{\bar{R}_1-\tau}{\tilde{c}}}^{\infty}u^{\frac{-1}{k+1}}u^{-n\alpha-1}du,\\
\end{split}
\end{equation}
where
\begin{equation}
\begin{split}
&\,\,\,\,\,\,\,\,\int_{\frac{\bar{R}_1-\tau}{\tilde{c}}}^{\infty}u^{\frac{-1}{k+1}}u^{-n\alpha-1}du\\
&=\int_{\frac{\bar{R}_1-\tau}{\tilde{c}}}^{\infty}u^{-n\alpha-1-\frac{1}{k+1}}du\\
&=\frac{1}{-n\alpha-\frac{1}{k+1}}((\bar{R}_1-\tau)/\tilde{c})^{-n\alpha-\frac{1}{k+1}}.\\
\end{split}
\end{equation}

\item \begin{equation}
\begin{split}
&\,\,\,\,\,\,\,\,\,\,\bar{T}_{k,2,N}\\
&=\int_{-\infty}^{\bar{R}_2}\bar{f}_k(y;\theta)\bar{p}_{2,\,N}((y-\tau)/\tilde{c})\frac{1}{\tilde{c}}dy\\
&=\frac{((\tau-\bar{R}_2)/\tilde{c})^{2}}{\pi}\sum_{n=1}^{N}(-1)^{n-1}\frac{\Gamma(n\alpha+1)}{n!}(\sin(\pi
n(1-\beta)/2\alpha))\\
&\,\,\,\,\,\,\,\cdot \int_{-\infty}^{\frac{\bar{R}_2-\tau}{\tilde{c}}}|u|^{\frac{-1}{k+1}}|u|^{-n\alpha-1}du,\\
\end{split}
\end{equation}
where
\begin{equation}
\begin{split}
&\,\,\,\,\,\,\,\,\int_{-\infty}^{\frac{\bar{R}_2-\tau}{\tilde{c}}}|u|^{\frac{-1}{k+1}}|u|^{-n\alpha-1}du\\
&=\int_{\frac{\tau-\bar{R}_2}{\tilde{c}}}^{\infty}u^{-(n+1)-\frac{1}{k+1}}du\\
&=\frac{1}{-n\alpha-\frac{1}{k+1}}((\tau-\bar{R}_2)/\tilde{c})^{-n\alpha-\frac{1}{k+1}}.\\
\end{split}
\end{equation}

\item \begin{equation}
\begin{split}
&\,\,\,\,\,\,\,\,\,\,\bar{T}_{k,3,N}\\
&=\int_{\bar{R}_1}^{\bar{R}_2}\bar{f}_k(y;\theta)p_{3,\,N}((y-\tau)/\tilde{c})\frac{1}{\tilde{c}}dy\\
&=\frac{1}{\pi}\sum_{n=1}^{N}(-1)^{n-1}\frac{\Gamma(n/\alpha+1)}{n!}(\sin(\pi
n(1+\beta)/2))\\
&\,\,\,\,\,\,\,\cdot \int_{\frac{\bar{R}_2-\tau}{\tilde{c}}}^{\frac{\bar{R}_1-\tau}{\tilde{c}}}u^{2(1-\frac{1}{2(k+1)})}u^{n-1}du,\\
\end{split}
\end{equation}
\end{enumerate}
where
\begin{equation}
\int_{\frac{\bar{R}_2-\tau}{\tilde{c}}}^{\frac{\bar{R}_1-\tau}{\tilde{c}}}u^{2(1-\frac{1}{2(k+1)})}u^{n-1}du=\\
\begin{cases}
\bar{A}_1+\bar{A}_2, &\text{ if $n$ is odd;}\\
\bar{A}_1-\bar{A}_2, &\text{ if $n$ is even,}\\
\end{cases}
\end{equation}
with
\begin{equation}
\begin{split}
\bar{A}_1&=\int_{0}^{\frac{\bar{R}_1-\tau}{\tilde{c}}}u^{2(1-\frac{1}{2(k+1)})}u^{n-1}du\\
   &=\int_{0}^{\frac{\bar{R}_1-\tau}{\tilde{c}}}u^{n-1+2(1-\frac{1}{2(k+1)})}du\\
   &=\frac{1}{n+2(1-\frac{1}{2(k+1)})}((\bar{R}_1-\tau)/\tilde{c})^{n+2(1-\frac{1}{2(k+1)})},
\end{split}
\end{equation}
and
\begin{equation}
\begin{split}
\bar{A}_2&=\int_{0}^{\frac{\tau-\bar{R}_2}{\tilde{c}}}u^{2(1-\frac{1}{2(k+1)})}u^{n-1}du\\
   &=\int_{0}^{\frac{\tau-\bar{R}_2}{\tilde{c}}}u^{n-1+2(1-\frac{1}{2(k+1)})}du\\
   &=\frac{1}{n+2(1-\frac{1}{2(k+1)})}((\tau-\bar{R}_2)/\tilde{c})^{n+2(1-\frac{1}{2(k+1)})}.
\end{split}
\end{equation}

Now we have the expressions of $\bar{T}_{k,1,N}$,
$\bar{T}_{k,2,N}$, and $\bar{T}_{k,3,N}$ without integrals inside.
We can substitute the sum of $\bar{T}_{k,1,N}$, $\bar{T}_{k,2,N}$
and $\bar{T}_{k,3,N}$ for $E(\bar{f}_k(x;\theta)),\,\,k=1,...,m$
in the expression of $\bar{V}$ of (3.3) and then do optimization
to get the estimate $\hat{\theta}$ for $\theta_0$.

\bigskip
\bigskip

{\bf Acknowledgement}: I am very grateful to Professor Anil K.
Bera for his encouragement and suggestions.

\vspace{.5in}
\begin{singlespace}

\end{singlespace}

\end{doublespace}

\end{document}